\documentclass[reqno,10pt]{amsart}

\usepackage{amsmath}
\usepackage{amssymb}
\usepackage{amsthm}
\usepackage{mathrsfs}
\usepackage{parskip}
\usepackage{geometry}
\usepackage{graphicx}
\usepackage{color}
\usepackage{graphics}
\usepackage[latin1]{inputenc}
\usepackage{verbatim}
\usepackage{tikz}
\usepackage{multirow}
\usepackage{expdlist}
\usepackage{stmaryrd}

\numberwithin{equation}{section}

\newtheorem{theorem}{Theorem}
\newtheorem{lemma}{Lemma}

\newtheorem{conjecture}{Conjecture}

\newcommand{\SC}{\operatorname{sc}}
\newcommand{\gn}{\mathcal{G}_n}
\newcommand{\leqRP}{\rotatebox{45}{$\leq$}}
\newcommand{\leqRN}{\rotatebox{-45}{$\leq$}}
\newcommand{\geqRP}{\rotatebox{45}{$\geq$}}
\newcommand{\geqRN}{\rotatebox{-45}{$\geq$}}

\newcommand{\lRP}{\rotatebox{45}{$<$}}
\newcommand{\lRN}{\rotatebox{-45}{$<$}}
\newcommand{\gRP}{\rotatebox{45}{$>$}}
\newcommand{\gRN}{\rotatebox{-45}{$>$}}
\newcommand{\eqRP}{\rotatebox{45}{$=$}}
\newcommand{\eqRN}{\rotatebox{-45}{$=$}}
\newcommand{\br}{(k_1,\ldots,k_n)}
\newcommand{\id}{\mathrm{id}}
\newcommand{\sts}{\large} 
\newcommand{\lbb}{\textbf{(}}
\newcommand{\rbb}{\textbf{)}}

\begin{document}

\author{Lukas Riegler}
\title[Generalized Monotone Triangles]{Generalized Monotone Triangles: an
extended Combinatorial Reciprocity Theorem} 
\thanks{Supported by the Austrian Science Foundation FWF, START grant Y463.}
\thanks{Fakultät für Mathematik, Universität Wien, Nordbergstraße 15, A-1090
Wien, Austria}
\thanks{E-Mail:lukas.riegler@univie.ac.at}

\begin{abstract}
In a recent work, the combinatorial interpretation of the polynomial $\alpha(n;
k_1,k_2,\ldots,k_n)$ counting the number of Monotone Triangles with bottom row
$k_1 < k_2 < \cdots < k_n$ was extended to weakly decreasing sequences $k_1
\geq k_2 \geq \cdots \geq k_n$. In this case the evaluation of the polynomial
is equal to a signed enumeration of objects called Decreasing Monotone
Triangles. In this paper we define Generalized Monotone Triangles -- a joint
generalization of both ordinary Monotone Triangles and Decreasing Monotone
Triangles. As main result of the paper we prove that the evaluation of
$\alpha(n; k_1,k_2,\ldots,k_n)$ at arbitrary $(k_1,k_2,\ldots,k_n) \in
\mathbb{Z}^n$ is a signed enumeration of Generalized Monotone Triangles with
bottom row $(k_1,k_2,\ldots,k_n)$. Computational experiments indicate that
certain evaluations of the polynomial at integral sequences yield well-known
round numbers related to Alternating Sign Matrices. The main result provides a combinatorial
interpretation of the conjectured identities and could turn out useful in giving
a bijective proof.
\end{abstract}

\maketitle

\section{\sts Introduction}
\begin{figure}[ht]
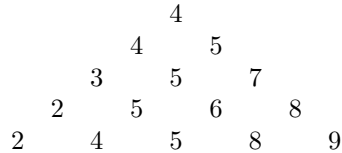

\begin{center}
$
\begin{array}{ccccccccc}
&&&& 4 \\
&&& 4 && 5  \\
&& 3 && 5 && 7 \\
& 2 && 5 && 6 && 8 \\
2 && 4 && 5 && 8 && 9 \\
\end{array}
$
\end{center}
\caption{One of the $16939$ Monotone Triangles with bottom row
$(2,4,5,8,9)$. \label{mt-example}}
\end{figure}

A \emph{Monotone Triangle} of size $n$ is a triangular
array of integers $(a_{i,j})_{1 \leq j \leq i \leq n}$
\[
\begin{array}{ccccccc}
&&& a_{1,1} \\
&& a_{2,1} && a_{2,2} \\
& \rotatebox{75}{$\ddots$} &&&& \ddots \\
a_{n,1} &&\cdots&  &\cdots&& a_{n,n}
\end{array}
\]
with strict increase along rows and weak increase along North-East- and
South-East-diagonals, i.e. $a_{i,j} < a_{i,j+1}$, $a_{i+1,j} \leq a_{i,j} \leq a_{i+1,j+1}$. An example of a Monotone Triangle of size $5$ is given in
Fig.\ref{mt-example}.

For each $n \geq 1$, there exists a unique polynomial
$\alpha(n;k_1,k_2,\ldots,k_n)$ of degree $n-1$ in each of the $n$ variables such that the evaluation of this polynomial at strictly increasing
sequences $k_1 < k_2 < \cdots < k_n$ is equal to the number of Monotone
Triangles with prescribed bottom row $(k_1,k_2,\ldots,k_n)$ -- for example
$\alpha(5;2,4,5,8,9) = 16939$.
This result was derived in \cite{FischerNumberOfMT}, where the polynomials are
given explicitly in terms of an operator formula.

In \cite{FischerRieglerDMT} we studied the evaluation of
$\alpha(n;k_1,\ldots,k_n)$ at weakly decreasing sequences $k_1 \geq k_2 \geq
\cdots \geq k_n$. It turned out that the evaluation can be interpreted as signed
enumeration of the following combinatorial objects:

A \emph{Decreasing Monotone Triangle} (DMT) of
size $n$ is a triangular array of integers $(a_{i,j})_{1 \leq j \leq i \leq n}$ having the following properties:
\begin{itemize}
\item\label{dmtDef1} The entries along North-East- and South-East-diagonals are
weakly decreasing.
\item\label{dmtDef2} Each integer appears at most twice in a row.
\item\label{dmtDef3} Two consecutive rows do not contain the same integer
exactly once.
\end{itemize}
One of the motivations for considering evaluations of
$\alpha(n;k_1,\ldots,k_n)$ at non-increasing $(k_1,\ldots,k_n)\in\mathbb{Z}^n$
stems from the connection to Alternating Sign Matrices. An \emph{Alternating Sign
Matrix} (ASM) of size $n$ is a $n \times n$-matrix with entries in $\{0,1,-1\}$
such that in each row and column the non-zero entries alternate in sign and sum up
to $1$. It is well-known that the set of ASMs is in bijection with the set of
Monotone Triangles with bottom row $(1,2,\ldots,n)$. Counting the number of ASMs
of size $n$ had been an open problem for more than a decade until the first two
independent proofs were given by D. Zeilberger (\cite{ZeilbergerASMProof}) and
G. Kuperberg (\cite{KuperbergASMProof}) in $1996$ (see \cite{Bressoud} for more
details). The Refined ASM Theorem -- i.e. the refined enumeration with respect
to the unique $1$ in the first row -- was reproven by I. Fischer in $2007$
(\cite{FischerNewRefProof}). The identity
\begin{equation}
\label{cylicEq}
\alpha(n;k_1,\ldots,k_n) = (-1)^{n-1} \alpha(n;k_2,\ldots,k_n,k_1-n)
\end{equation}
plays one of the key roles in this algebraic proof. A bijective proof
of \eqref{cylicEq} could give more combinatorial insight to the theorem.
However, note that if $k_1 < k_2 < \cdots < k_n$, then $k_n > k_1 - n$, i.e. \eqref{cylicEq} can per
se only be understood as identity satisfied by the polynomial.

The objective of this paper is to give an interpretation to the evaluation of
$\alpha(n;k_1,\ldots,k_n)$ at arbitrary $(k_1,\ldots,k_n) \in \mathbb{Z}^n$. For
this, we define triangular arrays of integers which locally combine the
restrictions of ordinary Monotone Triangles and Decreasing Monotone Triangles:

A \emph{Generalized Monotone Triangle} (GMT) is a
triangular array $(a_{i,j})_{1\leq j \leq i \leq n}$ of integers satisfying the
following conditions:
\begin{enumerate}
  \item\label{gmtDef1} Each entry is weakly bounded by its SW- and SE-neighbour,
  i.e.
  \[
  \min\{a_{i+1,j},a_{i+1,j+1}\} \leq a_{i,j} \leq \max\{a_{i+1,j},a_{i+1,j+1}\}.
  \]
  \item\label{gmtDef2} If three consecutive entries in a row are weakly
  increasing, then their two interlaced neighbours in the row above are strictly increasing, i.e.
  \[
  a_{i+1,j} \leq a_{i+1,j+1} \leq a_{i+1,j+2} \rightarrow a_{i,j} < a_{i,j+1}.
  \]
  \item\label{gmtDef3} If two consecutive entries in a row are strictly
  decreasing and their interlaced neighbour in the row above is equal to its
  SW-/SE-neighbour, then the interlaced neighbour has a left/right neighbour and is equal to it, i.e.
  \[
  a_{i,j} = a_{i+1,j} > a_{i+1,j+1} \rightarrow a_{i,j-1} = a_{i,j},
  \]
  \[
  a_{i+1,j} > a_{i+1,j+1} = a_{i,j} \rightarrow a_{i,j+1} = a_{i,j}.
  \]
\end{enumerate}

\begin{figure}[ht]
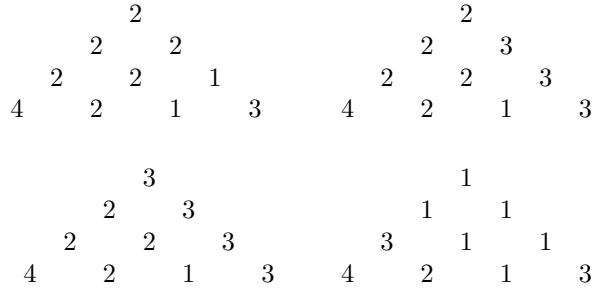

\begin{center}
$
\begin{array}{cc}
\begin{array}{ccccccccc}
&&& 2 \\
&& 2 && 2 \\
& 2 && 2 && 1 \\
4 && 2 && 1 && 3 \\ 
\end{array}
\quad &
\begin{array}{ccccccccc}
&&& 2 \\
&& 2 && 3 \\
& 2 && 2 && 3 \\
4 && 2 && 1 && 3 \\ 
\end{array}
\vspace{0.5cm}
\\
\begin{array}{ccccccccc}
&&& 3 \\
&& 2 && 3 \\
& 2 && 2 && 3 \\
4 && 2 && 1 && 3 \\ 
\end{array}
&
\begin{array}{ccccccccc}
&&& 1 \\
&& 1 && 1 \\
& 3 && 1 && 1 \\
4 && 2 && 1 && 3 \\ 
\end{array}
\end{array}
$
\end{center}
\caption{\label{gmtExample}The four GMTs with bottom row $(4,2,1,3)$.}
\end{figure}

By way of illustration, let us find all GMTs with bottom row $(4,2,1,3)$: First,
construct all possible penultimate rows $(l_1,l_2,l_3)$. Condition
\eqref{gmtDef1} implies that $l_1 \in \{2,3,4\}$, Condition \eqref{gmtDef3}
further restricts it to $l_1 \in \{2,3\}$. If on the one hand $l_1=2$, then
Condition \eqref{gmtDef3} forces $l_2=2$. The right-most entry $l_3$ is bounded
by $1$ and $3$, but actually $l_1=l_2=l_3=2$ would violate Condition
\eqref{gmtDef2}, so $l_3 \in \{1,3\}$. If on the other hand $l_1=3$, then
Condition \eqref{gmtDef3} implies that $l_2=l_3=1$. Continuing in the same way
with all penultimate rows yields the four GMTs depicted in Figure \ref{gmtExample}.

For $k_1 < k_2 < \cdots < k_n$, the set of GMTs with bottom row
$(k_1,\ldots,k_n)$ is equal to the set of Monotone Triangles with this bottom
row: Every GMT with strictly increasing bottom row is by
conditions \eqref{gmtDef1} and \eqref{gmtDef2} a Monotone Triangle. Conversely, the weak increase
along NE- and SE-diagonals of Monotone Triangles implies condition \eqref{gmtDef1} of
GMTs, the strict increase condition \eqref{gmtDef2}, and the premise of \eqref{gmtDef3} can not hold.

For $k_1 \geq k_2 \geq \cdots \geq k_n$, the set of GMTs with bottom row
$(k_1,\ldots,k_n)$ is equal to the set of Decreasing Monotone Triangles with
this bottom row: The NE- and SE-diagonals of every GMT with weakly decreasing
bottom row are by condition \eqref{gmtDef1} weakly decreasing. This also implies a weak
decrease along rows, and thus three consecutive equal entries in a row would contradict
condition \eqref{gmtDef2}. Furthermore, two consecutive rows containing an integer exactly
once would contradict condition \eqref{gmtDef3}. Conversely, the weak decrease of DMTs
along NE- and SE-diagonals implies condition \eqref{gmtDef1} and weak decrease along
rows. Thus, the premise of \eqref{gmtDef2} can only hold if three consecutive entries
coincide, which is not admissible in DMTs. Finally, condition \eqref{gmtDef3} follows from
the weak decrease along rows together with the condition that two consecutive
rows do not contain the same entry exactly once.

Therefore, Generalized Monotone Triangles are indeed a joint generalization of
ordinary Monotone Triangles and Decreasing Monotone Triangles. The main result
of the paper is that the evaluation $\alpha(n;k_1,\ldots,k_n)$ is a signed
enumeration of the GMTs with bottom row $(k_1,k_2,\ldots,k_n)$. The sign of a
GMT is determined by the following two statistics:

\begin{enumerate}
  \item An entry $a_{i,j}$ is called \emph{newcomer} if $a_{i+1,j} > a_{i,j} >
a_{i+1,j+1}$. 
\item A pair $(x,x)$ of two consecutive equal entries in a row is
called \emph{sign-changing}, if their interlaced neighbour in the row below is also equal to $x$.
\end{enumerate}

Let $\mathcal{G}_n(k_1,k_2,\ldots,k_n)$ denote the set of GMTs with bottom row
$(k_1,k_2,\ldots,k_n)$.
\begin{theorem}
\label{gmtTheorem}
Let $n \geq 1$ and $(k_1,k_2,\ldots,k_n)\in\mathbb{Z}^n$. Then
\[
\alpha(n;k_1,\ldots,k_n) = \sum_{A \in \mathcal{G}_n(k_1,\ldots,k_n)}
(-1)^{\SC(A)},
\]
where $\SC(A)$ is the total number of newcomers and sign-changing pairs in $A$. 
\end{theorem}

Applying Theorem \ref{gmtTheorem} to our example in Figure \ref{gmtExample}
yields $\alpha(4;4,2,1,3)=-2$.

Theorem \ref{gmtTheorem} is known to be true for strictly increasing
sequences $k_1 < k_2 < \cdots < k_n$, as in this case the set $\mathcal{G}_n(k_1,\ldots,k_n)$ is equal to the
set of Monotone Triangles with bottom row $(k_1,k_2,\ldots,k_n)$ and $\SC(A)=0$
for every Monotone Triangle.

Lemma $3$ of \cite{FischerRieglerDMT} implies the
correctness of Theorem \ref{gmtTheorem} for weakly decreasing bottom rows: In
this case $\mathcal{G}_n(k_1,\ldots,k_n)$ is equal to the set of DMTs with
bottom row $(k_1,\ldots,k_n)$ and the $\SC$-functions coincide. K. Jochemko and
R. Sanyal recently gave a proof of the theorem in this case from a geometric
point of view (\cite{JochemkoSanyalGeomProof}).

%
%

In Section \ref{sectionProof1} we give a straight-forward proof of Theorem
\ref{gmtTheorem} using a recursion satisfied by $\alpha(n;k_1,\ldots,k_n)$. In
Section \ref{sectionProof2} a connection with a known generalization (\cite{FischerSeqLabTrees}) is
established, which enables us to give a shorter, more subtle proof of Theorem
\ref{gmtTheorem}. Apart from being a joint generalization of Monotone Triangles
and DMTs, this generalization is more reduced in the sense that fewer
cancellations occur in the signed enumerations than in previously known
generalizations. In Section $4$ we apply the theorem to give a combinatorial
proof of an identity satisfied by $\alpha(n;k_1,\ldots,k_n)$ and provide a
collection of open problems.

\section{\sts Proof of Theorem \ref{gmtTheorem}}
\label{sectionProof1}

The number of Monotone Triangles with bottom row
$(k_1,\ldots,k_n)$ can be counted recursively by determining all admissible
penultimate rows $(l_1,\ldots,l_{n-1})$ and summing over the number of Monotone
Triangles with these bottom rows. The polynomial
$\alpha(n; k_1,\ldots,k_n)$ hence satisfies
\begin{equation}
\label{MTRec}
\alpha(n;k_1,\ldots,k_n)= \sum_{\substack{(l_1,\ldots,l_{n-1})\in
\mathbb{Z}^{n-1}, \\
k_1 \leq l_1 \leq k_2 \leq l_2 \leq \cdots \leq k_{n-1} \leq l_{n-1} \leq k_n,
\\ l_i < l_{i+1}}} \alpha(n-1; l_1,\ldots,l_{n-1})
\end{equation}
for all $k_1 < k_2 < \cdots < k_n$, $k_i \in \mathbb{Z}$. In fact
(\cite{FischerNumberOfMT}), one can define a summation operator
$\sum\limits_{(l_1,\ldots,l_{n-1})}^{(k_1,\ldots,k_n)}$ for arbitrary $(k_1,\ldots,k_n) \in \mathbb{Z}^n$ such that 
\begin{equation}
\label{alphaRec}
\alpha(n;k_1,\ldots,k_n) = \sum_{(l_1,\ldots,l_{n-1})}^{\br} \alpha(n-1;l_1,\ldots,l_{n-1})
\end{equation}
holds. The summation operator is defined recursively for arbitrary
$(k_1,\ldots,k_n) \in \mathbb{Z}^n$:
\begin{align}
\label{sumOpRec}
\sum_{(l_1,\ldots,l_{n-1})}^{\br} A(l_1,\ldots,l_{n-1}) := 
&\sum_{(l_1,\ldots,l_{n-2})}^{(k_1,\ldots,k_{n-1})}
\sum_{l_{n-1} = k_{n-1}+1}^{k_n} A(l_1,\ldots,l_{n-2},l_{n-1}) \\\notag &+
\sum_{(l_1,\ldots,l_{n-2})}^{(k_1,\ldots,k_{n-2},k_{n-1}-1)}
A(l_1,\ldots,l_{n-2},k_{n-1}), \quad n \geq 2,
\end{align}
with $\sum\limits_{()}^{(k_1)}:=\id$ and the extended definition of simple
sums
\begin{equation}
\label{sumExtDef}
\sum_{i=a}^b f(i) := \begin{cases} 0, & \quad b=a-1, \\ -\sum\limits_{i =
b+1}^{a-1} f(i), & \quad b+1 \leq a-1. \end{cases}
\end{equation}

Using induction and \eqref{MTRec}, it is clear that \eqref{alphaRec} holds
for increasing sequences $k_1 < k_2 < \cdots < k_n$. To prove it for arbitrary
$(k_1,\ldots,k_n)\in\mathbb{Z}^n$, let us first note that applying
$\sum\limits_{(l_1,\ldots,l_{n-1})}^{(k_1,\ldots,k_n)}$ to a polynomial in
$(l_1,\ldots,l_{n-1})$ yields a polynomial in $(k_1,\ldots,k_n)$: In the base
case $n=2$, write the polynomial $p(l_1)$ in terms of the binomial basis
$p(l_1)=\sum_{i=0}^{n-1} c_i \binom{l_1}{i}$. The polynomial
$q(x):=\sum_{i=1}^{n} c_{i-1} \binom{x}{i}$ then satisfies $q(x+1)-q(x)=p(x)$.
For integers $a \leq b$, it follows that $\sum_{l_1=a}^b p(l_1) = q(b+1)-q(a)$,
but this is by definition \eqref{sumExtDef} true for arbitrary $a,b \in
\mathbb{Z}$. The inductive step is immediate using \eqref{sumOpRec}. Thus, we
know that the right-hand side of \eqref{alphaRec} is a polynomial in
$(k_1,\ldots,k_n)$ coinciding with the polynomial on the left-hand side
whenever $k_1 < k_2 < \cdots < k_n$. Since a polynomial in $n$ variables is 
uniquely determined by these values, it follows that \eqref{alphaRec}
indeed holds.
The same is true for the alternative recursive description
\begin{align}
\label{sumOpRec2}
\sum_{(l_1,\ldots,l_{n-1})}^{\br} A(l_1,\ldots,l_{n-1}) =
& \sum_{(l_1,\ldots,l_{n-2})}^{(k_1,\ldots,k_{n-1})} \sum_{l_{n-1} =
k_{n-1}}^{k_n} A(l_1,\ldots,l_{n-2},l_{n-1}) \\\notag &-
\sum_{(l_1,\ldots,l_{n-3})}^{(k_1,\ldots,k_{n-2})}
A(l_1,\ldots,l_{n-3},k_{n-1},k_{n-1}), \quad n \geq 3.
\end{align}

The following Lemma establishes a connection between the summation operator and
GMTs, which then gives us the means to prove Theorem \ref{gmtTheorem}
inductively.

\begin{lemma}
\label{sumOpGMT}
Let $\mathcal{P}(k_1,\ldots,k_n)$ denote the set of $(n-1)$-st rows of elements
in $\gn(k_1,k_2,\ldots,k_n)$. Then every function $A(l_1,\ldots,l_{n-1})$
satisfies
\[
\sum_{(l_1,\ldots,l_{n-1})}^{(k_1,\ldots,k_n)} A(l_1,\ldots,l_{n-1}) =
\sum_{(l_1,\ldots,l_{n-1}) \in \mathcal{P}(k_1,\ldots,k_n)}
(-1)^{\SC(\mathbf{k};\mathbf{l})} A(l_1,\ldots,l_{n-1}), \quad n \geq 2,
\]
where $\SC(\mathbf{k};\mathbf{l}) :=
\SC(k_1,\ldots,k_n;l_1,\ldots,l_{n-1})$ is the total number of newcomers
and sign-changing pairs in $(l_1,\ldots,l_{n-1})$.
\end{lemma}
Before proving the Lemma, let us first give a remark, which is solely based on
the definition of GMTs. In general, the set of admissible values for an entry
$l_i$ depends on its neighbours $l_{i-1}$ and $l_{i+1}$ as well as the four
adjacent entries $k_{i-1}$, $k_i$, $k_{i+1}$ and $k_{i+2}$ in the row below
-- ordered
\begin{center}
$
\begin{array}{cccccccccccc}
& l_{i-1} && l_{i} && l_{i+1} \\
k_{i-1} && k_{i} && k_{i+1} && k_{i+2} \\ 
\end{array}
$
\end{center}
-- in the following way:
If $k_{i-1} > l_{i-1}=k_i$, then the only admissible value is $l_i = k_i$.
Symmetrically, if $k_{i+1} = l_{i+1} > k_{i+2}$, then $l_i =k_{i+1}$. Otherwise,
$l_i$ can take any value strictly between $k_i$ and $k_{i+1}$. To determine,
whether $l_i = k_i$ is allowed, check whether $k_i > k_{i+1}$, or $k_{i-1} > k_i
\leq k_{i+1}$ or $k_{i-1} \leq k_i \leq k_{i+1}$. If $k_{i} > k_{i+1}$, then
$l_i=k_i$ is admissible, if and only if $l_{i-1}=k_i$. If $k_{i-1} > k_i \leq k_{i+1}$,
then $l_i=k_i$ is admissible. If $k_{i-1} \leq k_i \leq k_{i+1}$, then
$l_i=k_i$ is admissible, if and only if $l_{i-1} < k_i$. Determining whether $l_i =
k_{i+1}$ is admissible works symmetrically.

\begin{proof}
%

If $n=2$, then the result is immediate using \eqref{sumExtDef}. For $n=3$, let
us check the case $k_1 \geq k_2$, $k_2 < k_3$ (the other cases can be shown
in the same way):
\begin{multline*}
\sum_{(l_1,l_2)}^{(k_1,k_2,k_3)} A(l_1,l_2) \stackrel{\eqref{sumOpRec2}}{=} 
\sum_{(l_1)}^{(k_1,k_2)} \sum_{l_2=k_2}^{k_3} A(l_1,l_2)-\sum_{()}^{(k_1)} A(k_2,k_2) \\
\stackrel{\eqref{sumOpRec}}{=}
 \sum_{()}^{(k_1)}\sum_{l_1=k_1+1}^{k_2}\sum_{l_2=k_2}^{k_3} A(l_1,l_2)+\sum_{()}^{(k_1-1)}\sum_{l_2=k_2}^{k_3} A(k_1,l_2)-A(k_2,k_2) \\
=
\sum_{l_1=k_1}^{k_2}\sum_{l_2=k_2}^{k_3}A(l_1,l_2)-A(k_2,k_2).
\end{multline*}
If $k_1=k_2$, then $\mathcal{P}(k_1,k_2,k_3)=\{(k_1,l_2):k_2<l_2 \leq k_3\}$
with no newcomers or sign-changing pairs. If $k_1 > k_2$, then
$\mathcal{P}(k_1,k_2,k_3)=\{(l_1,l_2):k_1>l_1>k_2, k_2 \leq l_2 \leq k_3 \}
\cup \{(k_2,k_2)\}$. The entry $l_1$ is either a newcomer or contained in
a sign-changing pair. The claimed equation follows from \eqref{sumExtDef}.

For $n \geq 4$, we have to distinguish between the cases $k_{n-1} \leq k_n$
(Case $1$) and $k_{n-1} > k_n$ (Case $2$). The remark preceding the proof also
suggests a different behaviour depending on whether $l_{n-1}$ -- the
rightmost entry of the penultimate row -- is equal to $k_{n-1}$ or
not. Indeed, this yields the sub-cases $1.1$, $1.2$ and $2.1$, $2.2$
respectively.

$\mathbf{\textbf{Case } 1 \; \lbb k_{n-1} \leq k_n \rbb :}$

Recursion \eqref{sumOpRec} of the summation operator and the
induction hypothesis yield
\begin{multline*}
\sum_{(l_1,\ldots,l_{n-1})}^{(k_1,\ldots,k_n)} A(l_1,\ldots,l_{n-1}) \\ =
\sum_{(l_1,\ldots,l_{n-2})}^{(k_1,\ldots,k_{n-1})}
\sum_{l_{n-1}=k_{n-1}+1}^{k_n} A(l_1,\ldots,l_{n-1}) +
\sum_{(l_1,\ldots,l_{n-2})}^{(k_1,\ldots,k_{n-2},k_{n-1}-1)}
A(l_1,\ldots,l_{n-2},k_{n-1}) \\
= \sum_{(l_1,\ldots,l_{n-2})\in \mathcal{P}(k_1,\ldots,k_{n-1})}
(-1)^{\SC(k_1,\ldots,k_{n-1};l_1,\ldots,l_{n-2})} \sum_{l_{n-1}=k_{n-1}+1}^{k_n}
A(l_1,\ldots,l_{n-1}) \\ +
\sum_{(l_1,\ldots,l_{n-2})\in \mathcal{P}(k_1,\ldots,k_{n-2},k_{n-1}-1)}
(-1)^{\SC(k_1,\ldots,k_{n-2},k_{n-1}-1;l_1,\ldots,l_{n-2})}
A(l_1,\ldots,l_{n-2},k_{n-1}).
\end{multline*}

To see that this is further equal to
\[
\sum_{(l_1,\ldots,l_{n-1}) \in \mathcal{P}(k_1,\ldots,k_n)}
(-1)^{\SC(\mathbf{k};\mathbf{l})} A(l_1,\ldots,l_{n-1}),
\]
let us show that for $k_{n-1} \leq k_n$
\begin{multline}
\label{penIso1}
\mathcal{P}(k_1,\ldots,k_n)  \\ 
= \mathcal{P}(k_1,\ldots,k_{n-1}) \times \{l_{n-1} \mid k_{n-1} < l_{n-1} \leq
k_n \} \; \cup \; \mathcal{P}(k_1,\ldots,k_{n-2},k_{n-1}-1) \times
\{k_{n-1}\} \quad
\end{multline}
holds, and that each fixed row causes the same total number of sign-changes on 
the left-hand side as on the right-hand side.

Since the first $n-2$ entries of the bottom row are identical on both
sides of \eqref{penIso1}, it suffices -- by the remark preceding the proof -- to
show that the restrictions imposed on $l_{n-1}$, $l_{n-2}$ and $l_{n-3}$ are the same on both
sides. For this, consider the entry $l_{n-1}$ and distinguish between $k_{n-1}
< l_{n-1} \leq k_n$ and $l_{n-1} = k_{n-1}$:

$\mathbf{\textbf{Case } 1.1 \; \lbb k_{n-1} < l_{n-1} \leq k_n \rbb :}$

If $k_{n-2} > k_{n-1}$, then $k_{n-2} \geq l_{n-2} >
k_{n-1}$ on both sides:
\begin{center}
$
\begin{array}{|c|}
\hline 
\text{Left-hand side of \eqref{penIso1}} \\ 
\hline
\begin{array}{cccccccccc}
&& l_{n-2} &&&&& l_{n-1} \\
& \geqRP && \gRN && \lRP &&& \leqRN \\
k_{n-2} && > && k_{n-1} &&& \leq && k_n \\ 
\end{array}
\\
\hline\hline
\text{Right-hand side of \eqref{penIso1}}\\
\hline
\begin{array}{ccccccccc}
&& l_{n-2} &&&&& l_{n-1} \\
& \geqRP && \gRN \\
k_{n-2} && > && k_{n-1} \\
\end{array}
\\
\hline
\end{array}
$
\end{center}
The restrictions for $l_{n-3}$ are the same on both sides. The entry $l_{n-1}$ does not contribute a sign-change,  and the entry $l_{n-2}$
is involved in a sign-change on both sides.

If $k_{n-2} \leq k_{n-1}$,
then $k_{n-2} \leq l_{n-2} \leq k_{n-1}$ on both sides:
\begin{center}
$
\begin{array}{|c|}
\hline 
\text{Left-hand side of \eqref{penIso1}} \\ 
\hline
\begin{array}{cccccccccc}
&& l_{n-2} &&&&& l_{n-1} \\
& \leqRP && \leqRN && \lRP &&& \leqRN \\
k_{n-2} && \leq && k_{n-1} &&& \leq && k_n \\ 
\end{array}
\\
\hline\hline
\text{Right-hand side of \eqref{penIso1}} \\
\hline
\begin{array}{ccccccccc}
&& l_{n-2} &&&&& l_{n-1} \\
& \leqRP && \leqRN \\
k_{n-2} && \leq && k_{n-1} \\
\end{array}
\\
\hline
\end{array}
$
\end{center}
The restrictions for $l_{n-3}$ are the same on both sides. The entry $l_{n-1}$ does not contribute a sign-change, and the entry $l_{n-2}$
is involved in a sign-change on the left-hand side, if and only if it is on the
right-hand side. 

It follows that $\SC(k_1,\ldots,k_{n-1};l_1,\ldots,l_{n-2}) =
\SC(k_1,\ldots,k_n;l_1,\ldots,l_{n-1})$.

$\mathbf{\textbf{Case } 1.2 \; \lbb l_{n-1} = k_{n-1} \rbb :}$

If $k_{n-2} = k_{n-1}$, then there is no row on the
left-hand side with $l_{n-1} = k_{n-1}$, and on the right-hand side this would
imply $l_{n-3} = l_{n-2} = l_{n-1} = k_{n-1}$:
\begin{center}
$
\begin{array}{|c|}
\hline 
\text{Left-hand side of \eqref{penIso1}} \\ 
\hline
\begin{array}{cccccccccc}
&& \lightning &&&&& k_{n-1} \\
&  &&  && \eqRP &&& \leqRN \\
k_{n-2} && = && k_{n-1} &&& \leq && k_n \\ 
\end{array}
\\
\hline\hline
\text{Right-hand side of \eqref{penIso1}} \\
\hline 
\begin{array}{cccccccccccccccc}
&& k_{n-1} &&&&& k_{n-1} &&&&& k_{n-1} \\
&  && \eqRN && \eqRP &&& \gRN \\
k_{n-3} &&  && k_{n-1} &&& > && k_{n-1}-1 \\ 
\end{array}
\\
\hline
\end{array}
$
\end{center}
But since a GMT can not contain three consecutive equal entries, such rows are
not contained on the right-hand side.

If $k_{n-2} \leq k_{n-1}-1$, then
$k_{n-2} \leq l_{n-2} < k_{n-1}$ on both sides:

\begin{center}
$
\begin{array}{|c|}
\hline 
\text{Left-hand side of \eqref{penIso1}} \\ 
\hline
\begin{array}{cccccccccc}
&& \l_{n-2} &&&&& k_{n-1} \\
& \leqRP && \lRN && \eqRP &&& \leqRN \\
k_{n-2} && < && k_{n-1} &&& \leq && k_n \\ 
\end{array}
\\
\hline\hline
\text{Right-hand side of \eqref{penIso1}} \\
\hline
\begin{array}{cccccccccc}
&&& l_{n-2} &&&&& k_{n-1} \\
& \leqRP &&& \leqRN \\
k_{n-2} &&& \leq && k_{n-1}-1 \\ 
\end{array}
\\
\hline
\end{array}
$
\end{center}
The restrictions for $l_{n-3}$ are the same on both sides. The entry $l_{n-1}$ does not contribute a sign-change, and the entry $l_{n-2}$
is involved in a sign-change on the left-hand side, if and only if it is on the
right-hand side.

If $k_{n-2} > k_{n-1}$, then $k_{n-2} \geq l_{n-2} \geq k_{n-1}$ on both sides:

\begin{center}
$
\begin{array}{|c|}
\hline 
\text{Left-hand side of \eqref{penIso1}} \\ 
\hline
\begin{array}{cccccccccc}
&& \l_{n-2} &&&&& k_{n-1} \\
& \geqRP && \geqRN && \eqRP &&& \leqRN \\
k_{n-2} && > && k_{n-1} &&& \leq && k_n \\ 
\end{array}
\\
\hline\hline
\text{Right-hand side of \eqref{penIso1}} \\
\hline
\begin{array}{cccccccccc}
&&& l_{n-2} &&&&& k_{n-1} \\
& \geqRP &&& \gRN \\
k_{n-2} &&& > && k_{n-1}-1 \\ 
\end{array}
\\
\hline
\end{array}
$
\end{center}
The restrictions for $l_{n-3}$ are the same on both sides. The entry $l_{n-2}$
is involved in a sign-change on both sides (note the special case
$l_{n-2}=k_{n-1}$, where $l_{n-2}$ is part of a sign-changing pair on the
left-hand side and a newcomer on the right-hand side).

It follows that $\SC(k_1,\ldots,k_{n-2},k_{n-1}-1;l_1,\ldots,l_{n-2}) =
\SC(k_1,\ldots,k_n;l_1,\ldots,l_{n-1})$.

$\mathbf{\textbf{Case } 2 \; \lbb k_{n-1} > k_n \rbb:}$

Recursion \eqref{sumOpRec2} of the summation
operator and the induction hypothesis yield
\begin{multline*}
\sum_{(l_1,\ldots,l_{n-1})}^{(k_1,\ldots,k_n)} A(l_1,\ldots,l_{n-1}) \\ =
-\sum_{(l_1,\ldots,l_{n-2})}^{(k_1,\ldots,k_{n-1})}
\sum_{l_{n-1}=k_n+1}^{k_{n-1}-1} A(l_1,\ldots,l_{n-1}) -
\sum_{(l_1,\ldots,l_{n-3})}^{(k_1,\ldots,k_{n-2})} A(l_1,\ldots,l_{n-3},k_{n-1},k_{n-1}) \\
=\sum_{(l_1,\ldots,l_{n-2}) \in \mathcal{P}(k_1,\ldots,k_{n-1})}
(-1)^{\SC(k_1,\ldots,k_{n-1};l_1,\ldots,l_{n-2})+1}
\sum_{l_{n-1}=k_n+1}^{k_{n-1}-1} A(l_1,\ldots,l_{n-1}) \\ +
\sum_{(l_1,\ldots,l_{n-3})\in \mathcal{P}(k_1,\ldots,k_{n-2})}
(-1)^{\SC(k_1,\ldots,k_{n-2};l_1,\ldots,l_{n-3})+1}
A(l_1,\ldots,l_{n-3},k_{n-1},k_{n-1}).
\end{multline*}

Similarly, let us show that for $k_{n-1}> k_n$
\begin{multline}
\label{penIso2}
\mathcal{P}(k_1,\ldots,k_n)  \\ 
= \mathcal{P}(k_1,\ldots,k_{n-1}) \times \{l_{n-1} \mid k_{n-1} > l_{n-1} >
k_{n} \} \; \cup \; \mathcal{P}(k_1,\ldots,k_{n-2}) \times \{(k_{n-1},k_{n-1})\} \quad
\end{multline}
holds. Again it suffices to show that
$l_{n-1}$, $l_{n-2}$ and $l_{n-3}$ have to satisfy the same restrictions on both
sides, and that corresponding rows contain the same number of
sign-changes. Since $k_{n-1} > k_n$, it follows that $k_{n-1}
\geq l_{n-1} > k_n$ on both sides. Let us distinguish between the cases $k_{n-1} > l_{n-1} > k_{n}$ and
$l_{n-1} = k_{n-1}$:

$\mathbf{\textbf{Case } 2.1 \; \lbb k_{n-1} > l_{n-1} > k_{n} \rbb :}$

If $k_{n-2} > k_{n-1}$, then $k_{n-2} \geq
l_{n-2} > k_{n-1}$ on both sides:

\begin{center}
$
\begin{array}{|c|}
\hline 
\text{Left-hand side of \eqref{penIso2}} \\ 
\hline
\begin{array}{cccccccccc}
&& \l_{n-2} &&&&& l_{n-1} \\
& \geqRP && \gRN && \gRP &&& \gRN \\
k_{n-2} && > && k_{n-1} &&& > && k_n \\ 
\end{array}
\\
\hline\hline
\text{Right-hand side of \eqref{penIso2}} \\
\hline
\begin{array}{cccccccccc}
&&& l_{n-2} &&&&& l_{n-1} \\
& \geqRP &&& \gRN \\
k_{n-2} &&& > && k_{n-1} \\ 
\end{array}
\\
\hline
\end{array}
$
\end{center}

The restrictions for $l_{n-3}$ are the same on both sides. The entries
$l_{n-1}$ and $l_{n-2}$ both contribute a sign-change.

If $k_{n-2} \leq k_{n-1}$, then $k_{n-2} \leq l_{n-2} \leq k_{n-1}$ on both
sides:

\begin{center}
$
\begin{array}{|c|}
\hline 
\text{Left-hand side of \eqref{penIso2}}\\ 
\hline
\begin{array}{cccccccccc}
&& \l_{n-2} &&&&& l_{n-1} \\
& \leqRP && \leqRN && \gRP &&& \gRN \\
k_{n-2} && \leq && k_{n-1} &&& > && k_n \\ 
\end{array}
\\
\hline\hline
\text{Right-hand side of \eqref{penIso2}} \\
\hline
\begin{array}{cccccccccc}
&&& l_{n-2} &&&&& l_{n-1} \\
& \leqRP &&& \leqRN \\
k_{n-2} &&& \leq && k_{n-1} \\ 
\end{array}
\\
\hline
\end{array}
$
\end{center}

The restrictions for $l_{n-3}$ are the same on both sides. The entry $l_{n-1}$
contributes a sign-change, and the entry $l_{n-2}$ is involved in a sign-change
on the left-hand side, if and only if it is on the right-hand side.

It follows that
$\SC(k_1,\ldots,k_{n-1};l_1,\ldots,l_{n-2})+1=\SC(k_1,\ldots,k_n;l_1,\ldots,l_{n-1})$.

$\mathbf{\textbf{Case } 2.2 \; \lbb l_{n-1} = k_{n-1} \rbb :}$

Since $k_{n-1} > k_n$ and $l_{n-1}=k_{n-1}$, we have $l_{n-2}=k_{n-1}$ on both
sides, whereby $(l_{n-2},l_{n-1})$ is a sign-changing pair. It remains to
be shown that $l_{n-3}$ has the same restrictions on both sides.

If $k_{n-3} \leq k_{n-2}$, then $k_{n-3} \leq l_{n-3} \leq k_{n-2}$ on the
right-hand side:

\begin{center}
$
\begin{array}{|c|}
\hline 
\text{Left-hand side of \eqref{penIso2}} \\ 
\hline
\begin{array}{ccccccccccccc}
&& l_{n-3} && && k_{n-1} &&&& k_{n-1} \\
& \leqRP &&  &&  && \eqRN && \eqRP && \gRN \\
k_{n-3} && \leq && k_{n-2} &&  && k_{n-1} && > && k_n \\ 
\end{array}
\\
\hline\hline
\text{Right-hand side of \eqref{penIso2}}
\\
\hline
\begin{array}{ccccccccccccccc}
&&& l_{n-3} &&&&& k_{n-1} &&&&& k_{n-1} \\
& \leqRP &&& \leqRN \\
k_{n-3} &&& \leq && k_{n-2} \\ 
\end{array}
\\
\hline
\end{array}
$
\end{center}
On the left-hand side we also have $k_{n-3} \leq l_{n-3} \leq k_{n-2}$, unless
$k_{n-2} = k_{n-1}$. In this case $k_{n-3} \leq l_{n-3} < k_{n-2}$, but for
$l_{n-3} = k_{n-2} = k_{n-1} = l_{n-2} = l_{n-1}$ there are three consecutive
equal entries anyway. The entry $l_{n-3}$ is involved in a sign-change on the left-hand side, if and
only if it is on the right-hand side, and
$(l_{n-2},l_{n-1})$ is a sign-changing pair.

If $k_{n-3} > k_{n-2}$ (and $n > 4$), then $k_{n-3} \geq l_{n-3} > k_{n-2}$ on
the right-hand side: 
\begin{center}
$
\begin{array}{|c|}
\hline 
\text{Left-hand side of \eqref{penIso2}} \\ 
\hline
\begin{array}{ccccccccccccc}
&& l_{n-3} && && k_{n-1} &&&& k_{n-1} \\
& \geqRP && &&  && \eqRN && \eqRP && \gRN \\
k_{n-3} && > && k_{n-2} &&  && k_{n-1} && > && k_n \\ 
\end{array}
\\
\hline\hline
\text{Right-hand side of \eqref{penIso2}}
\\
\hline
\begin{array}{ccccccccccccccc}
&&& l_{n-3} &&&&& k_{n-1} &&&&& k_{n-1}\\
& \geqRP &&& \gRN \\
k_{n-3} &&& > && k_{n-2} \\ 
\end{array}
\\
\hline
\end{array}
$
\end{center}
Again, $l_{n-3}$ has the same restrictions on the left-hand side, unless
$k_{n-2} = k_{n-1}$. In this case $k_{n-3} \geq l_{n-3} \geq k_{n-2}$, whereby
$l_{n-3} = k_{n-2}$ implies that $l_{n-3} = k_{n-2} = k_{n-1} = l_{n-2} =
l_{n-1}$. If $n=4$, then the same holds with the difference that $k_{1} > l_{1}$ instead of $k_1 \geq
l_1$ on both sides. The entry $l_{n-3}$ is involved in a sign-change on both sides and
$(l_{n-2},l_{n-1})$ is a sign-changing pair. 

It follows that
$\SC(k_1,\ldots,k_{n-2};l_1,\ldots,l_{n-3})+1=\SC(k_1,\ldots,k_n;l_1,\ldots,l_{n-1})$.

\end{proof}


\begin{proof}[Proof (Theorem \ref{gmtTheorem})]
The result is immediate for $n=1$:
\[
\alpha(1;k_1) = 1 = \sum_{A \in \mathcal{G}_1(k_1)} (-1)^{\SC(A)}.
\]
For $n \geq 2$ apply \eqref{alphaRec}, Lemma \ref{sumOpGMT} and the induction
hypothesis:
\begin{multline*}
\alpha(n;k_1,\ldots,k_n) = \sum_{(l_1,\ldots,l_{n-1})}^{(k_1,\ldots,k_n)} \alpha(n-1;l_1,\ldots,l_{n-1}) \\
= \sum_{(l_1,\ldots,l_{n-1}) \in \mathcal{P}(k_1,\ldots,k_n)} (-1)^{\SC(\mathbf{k};\mathbf{l})} \alpha(n-1;l_1,\ldots,l_{n-1}) \\
= \sum_{(l_1,\ldots,l_{n-1}) \in \mathcal{P}(k_1,\ldots,k_n)}
(-1)^{\SC(\mathbf{k};\mathbf{l})} \sum_{A \in
\mathcal{G}_{n-1}(l_1,\ldots,l_{n-1})} (-1)^{\SC(A)} = \sum_{A \in \mathcal{G}_n(k_1,\ldots,k_n)} (-1)^{\SC(A)}. \\
\end{multline*}

\end{proof}

\section{\sts Connection with different generalization \& Alternative proof}
\label{sectionProof2}

In \cite{FischerSeqLabTrees} four different combinatorial extensions of
$\alpha(n;k_1,\ldots,k_n)$ to all $(k_1,\ldots,k_n)\in\mathbb{Z}^n$ are
described. The idea behind all of them is to write the sum in \eqref{MTRec} in
terms of simple summations, i.e. summations as defined in \eqref{sumExtDef}.
In the third extension this is based on the inclusion-exclusion principle: Let
$k_1 < k_2 < \cdots < k_n$ and
\begin{align*}
M&:=\{(l_1,\ldots,l_{n-1}) \in \mathbb{Z}^{n-1}\mid \forall j: k_j \leq l_j \leq
k_{j+1} \; \land \; l_{j}<l_{j+1} \},
\\
A&:=\{(l_1,\ldots,l_{n-1}) \in \mathbb{Z}^{n-1} \mid \forall j: k_j \leq l_j
\leq k_{j+1}\}, \\
A_i&:=\{(l_1,\ldots,l_{n-1}) \in \mathbb{Z}^{n-1} \mid \forall j: k_j \leq l_j
\leq k_{j+1} \;\land\; l_{i-1}=k_i=l_i\}, \quad i=2,\ldots,n-1.
\end{align*}
The strict increase implies that $A_i \cap A_{i+1} = \emptyset$, and thus we
have for any function $f(\mathbf{l}):=f(l_1,\ldots,l_{n-1})$ that
\begin{multline}
\label{thirdExtConnMT}
\sum_{\mathbf{l}\in M} f(\mathbf{l}) =\sum_{\mathbf{l}\in
A}f(\mathbf{l})-\sum_{i=2}^{n-1}\sum_{\mathbf{l}\in
A_i}f(\mathbf{l})+\sum_{\substack{2 \leq i_1 < i_2 \leq n-1 \; \\i_2 \neq
i_1+1}}\sum_{\mathbf{l}\in A_{i_1}\cap A_{i_2}}f(\mathbf{l}) \\-
\sum_{\substack{2 \leq i_1 < i_2 < i_3 \leq n-1 \; \\i_{j+1} \neq i_j+1}}\sum_{\mathbf{l}\in A_{i_1}\cap A_{i_2}\cap A_{i_3}}f(\mathbf{l}) \cdots,
\end{multline}
which can be written in terms of simple sums as
\begin{equation}
\label{thirdExt}
\sum_{p \geq 0} (-1)^p\sum_{\substack{2 \leq i_1 < i_2 < \cdots < i_p \leq
n-1 \; \\i_{j+1} \neq i_j +1}}
\sum_{l_1=k_1}^{k_2}\sum_{l_2=k_2}^{k_3}\cdots\sum_{l_{i_1}-1=k_{i_1}}^{k_{i_1}}\sum_{l_{i_1}=k_{i_1}}^{k_{i_1}}
\cdots \sum_{l_{i_p}-1=k_{i_p}}^{k_{i_p}}\sum_{l_{i_p}=k_{i_p}}^{k_{i_p}}
\cdots \sum_{l_{n-1}=k_{n-1}}^{k_n} f(\mathbf{l}).
\end{equation}
Using \eqref{sumExtDef}, we can interpret \eqref{thirdExt} for arbitrary
$(k_1,\ldots,k_n) \in \mathbb{Z}^n$. Let us show that 
\begin{align}
\label{thirdExtRec}
\alpha(n;k_1,\ldots,k_n) &= \sum_{p \geq 0} (-1)^p\sum_{\substack{2 \leq i_1 <
i_2 < \cdots < i_p \leq n-1 \; \\i_{j+1} \neq i_j +1}} \\ \notag
&
\sum_{l_1=k_1}^{k_2}\sum_{l_2=k_2}^{k_3}\cdots\sum_{l_{i_1}-1=k_{i_1}}^{k_{i_1}}\sum_{l_{i_1}=k_{i_1}}^{k_{i_1}} \cdots \sum_{l_{i_p}-1=k_{i_p}}^{k_{i_p}}\sum_{l_{i_p}=k_{i_p}}^{k_{i_p}}
\cdots \sum_{l_{n-1}=k_{n-1}}^{k_n} \alpha(n-1;l_1,\ldots,l_{n-1})
\end{align}
holds for $(k_1,\ldots,k_n)\in\mathbb{Z}^n$. The correctness for $k_1 < k_2 < \cdots <
k_n$ is ensured by \eqref{MTRec}, \eqref{thirdExtConnMT} and \eqref{thirdExt}.
To prove it for arbitrary $(k_1,\ldots,k_n) \in \mathbb{Z}^n$, it thus suffices
to show that \eqref{thirdExt} applied to a polynomial in $l_1,\ldots,l_{n-1}$
yields a polynomial in $k_1,\ldots,k_n$. But this follows from
\eqref{sumExtDef} in the exact same way as in the proof of \eqref{alphaRec}. 

As pointed out in \cite{FischerSeqLabTrees}, we can give
\eqref{thirdExtRec} a combinatorial meaning by interpreting
$\alpha(n;k_1,\ldots,k_n)$ as signed enumeration of the following combinatorial
objects: In a triangular array $(a_{i,j})_{1 \leq j \leq i \leq n}$ of integers, let us call the entries $a_{i-1,j-1}$ and $a_{i-1,j}$ the parents of $a_{i,j}$. Among the entries $(a_{i,j})_{1 < j < i \leq n}$, there may be \emph{special entries}.
Special entries in the same row must not be adjacent (choosing these special
entries corresponds to fixing the $i_l$'s in \eqref{thirdExtRec}). The
requirements for the entries are
\begin{enumerate}
  \item If $a_{i,j}$ is special, then $a_{i-1,j-1} = a_{i,j} = a_{i-1,j}$.
  \item If $a_{i,j}$ is not the parent of a special entry and $a_{i+1,j} \leq
  a_{i+1,j+1}$, then $a_{i+1,j} \leq a_{i,j} \leq a_{i+1,j+1}$.
  \item If $a_{i,j}$ is not the parent of a special entry and $a_{i+1,j} >
  a_{i+1,j+1}$, then $a_{i+1,j+1} > a_{i,j} > a_{i+1,j}$. In this case $a_{i,j}$
  is called \emph{inversion}.
\end{enumerate}

Let us denote by $\mathcal{T}_n(k_1,\ldots,k_n)$ the set of these objects with
bottom row $(a_{n,1},\ldots,a_{n,n})=(k_1,\ldots,k_n)$. For $A \in
\mathcal{T}_n(k_1,\ldots,k_n)$ let $s(A)$ be the total number of special entries
and inversions. Using induction and \eqref{thirdExtRec}, we thus have
\[
\alpha(n;k_1,\ldots,k_n)=\sum_{A \in \mathcal{T}_n(k_1,\ldots,k_n)} (-1)^{s(A)}.
\]
We can now eliminate those arrays $(a_{i,j})_{1 \leq j \leq i \leq n}$
violating the condition
\begin{equation}
\label{incCond}
a_{i,j-1} \leq a_{i,j} \leq a_{i,j+1} \; \rightarrow \; a_{i-1,j-1} <
a_{i-1,j}
\end{equation}
by using the following sign-reversing involution: find the minimal index $i$,
and under those the minimal index $j$ such that $a_{i,j-1} \leq a_{i,j} \leq
a_{i,j+1}$ and $a_{i-1,j-1} = a_{i,j} = a_{i-1,j}$. If $a_{i,j}$ is special,
then turn it non-special, and vice-versa. Note that the minimality of $i$
ensures that turning $a_{i,j}$ special is admissible: Suppose a neighbour of
$a_{i,j}$ is special, then the row above contains three consecutive equal
entries and thus an entry violating \eqref{incCond}. It follows that 
\[
\alpha(n;k_1,\ldots,k_n)=\sum_{\substack{A \in \mathcal{T}_n(k_1,\ldots,k_n)
\\a_{i,j-1} \leq a_{i,j} \leq a_{i,j+1} \rightarrow a_{i-1,j-1} <
a_{i-1,j}} } (-1)^{s(A)}.
\]
Note that in this reduced set an entry $a_{i,j}$ is special if and only if
$a_{i-1,j-1}=a_{i,j}=a_{i-1,j}$. Hence, the additional information of which
entries are special is not required anymore. Since special entries now
correspond to sign-changing pairs and inversions to newcomers, the only remaining part for
proving Theorem \ref{gmtTheorem} is to show that
\[
\mathcal{G}_n(k_1,\ldots,k_n) = \{A \in \mathcal{T}_n(k_1,\ldots,k_n): a_{i,j-1}
\leq a_{i,j} \leq a_{i,j+1} \; \rightarrow \; a_{i-1,j-1} < a_{i-1,j} \},
\] 
where an entry $a_{i,j}$ is special if and only if
$a_{i-1,j-1}=a_{i,j}=a_{i-1,j}$.

Let $A \in \mathcal{G}_n(k_1,\ldots,k_n)$. Then two adjacent special entries in
a row would imply three consecutive equal entries in a row, in contradiction to
Condition \eqref{gmtDef2} of GMTs. If $a_{i,j}$ is special, then $a_{i-1,j-1} = a_{i,j} =
a_{i-1,j}$ by definition. If $a_{i+1,j} \leq a_{i+1,j+1}$, then $a_{i+1,j} \leq
a_{i,j} \leq a_{i+1,j+1}$ by Condition \eqref{gmtDef1} of GMTs. If $a_{i+1,j} >
a_{i+1,j+1}$, then $a_{i+1,j} \geq a_{i,j} \geq
a_{i+1,j+1}$ by Condition \eqref{gmtDef1} of GMTs, and if $a_{i+1,j}$ and
$a_{i+1,j+1}$ are neither special, Condition \eqref{gmtDef3} of GMTs implies that
$a_{i+1,j} > a_{i,j} > a_{i+1,j+1}$. We thus have $A \in \mathcal{T}_n(k_1,\ldots,k_n)$, and the
additional property is exactly Condition \eqref{gmtDef2} of GMTs.

Let $A \in \mathcal{T}_n(k_1,\ldots,k_n)$ such that $a_{i,j-1}
\leq a_{i,j} \leq a_{i,j+1}$ implies $a_{i-1,j-1} < a_{i-1,j}$. Conditions \eqref{gmtDef1}
and \eqref{gmtDef2} of GMTs are then trivially satisfied. If $a_{i,j} = a_{i+1,j} >
a_{i+1,j+1}$, then by Condition \eqref{gmtDef3} the entry $a_{i,j}$ has to be parent of a
special entry, and thus $a_{i,j} = a_{i+1,j} = a_{i,j-1}$. The second part of
Condition \eqref{gmtDef3} of GMTs is symmetric, and therefore $A \in
\mathcal{G}_n(k_1,\ldots,k_n)$.

This concludes the less straight-forward, yet much shorter proof of Theorem
\ref{gmtTheorem}.

\section{\sts Applications \& Open Problems}

With this generalization at hand, we can try to give a combinatorial
interpretation to identities satisfied by $\alpha(n;k_1,\ldots,k_n)$. By way of
illustration, take the identity
\begin{align}
\label{sasId1}
&\alpha(n;k_1,\ldots,k_{i-1},k_i,k_i +
1,k_{i+2},\ldots,k_n)\\ \notag
&=\alpha(n;k_1,\ldots,k_{i-1},k_i,k_i,k_{i+2},\ldots,k_n)+\alpha(n;k_1,\ldots,k_{i-1},k_i
+1,k_i + 1,k_{i+2},\ldots,k_n).
\end{align}
A combinatorial proof of this identity in the case that $k_1 < k_2 < \cdots <
k_i$ and $k_i+1 < k_{i+2} < \cdots < k_n$ was given in
\cite{FischerSeqLabTrees}. Using Theorem \ref{gmtTheorem}, we can now give a
combinatorial proof for arbitrary $(k_1,\ldots,k_n)\in\mathbb{Z}^n$ by showing
that there exists a sign-preserving bijection
\begin{multline*}
\mathcal{G}_n(k_1,\ldots,k_{i-1},k_i,k_i +
1,k_{i+2},\ldots,k_n) \\ \leftrightarrow
\mathcal{G}_n(k_1,\ldots,k_{i-1},k_i,k_i,k_{i+2},\ldots,k_n) \;\dot{\cup}\;
\mathcal{G}_n(k_1,\ldots,k_{i-1},k_i+1,k_i + 1,k_{i+2},\ldots,k_n).
\end{multline*}
If $\mathcal{P}(k_1,\ldots,k_n)$ denotes the set of penultimate rows of GMTs
with bottom row $(k_1,\ldots,k_n)$, it suffices to show that
\begin{multline}
\label{oneStepIso}
\mathcal{P}(k_1,\ldots,k_{i-1},k_i,k_i+1,k_{i+2},\ldots,k_n) \\
= \mathcal{P}(k_1,\ldots,k_{i-1},k_i,k_i,k_{i+2},\ldots,k_n) \;
\dot{\cup} \; \mathcal{P}(k_1,\ldots,k_{i-1},k_i+1,k_i+1,k_{i+2},\ldots,k_n),
\end{multline}
where each fixed row has the same total number of sign-changes on both sides.

Each $(l_1,\ldots,l_{n-1})\in
\mathcal{P}(k_1,\ldots,k_{i-1},k_i,k_i+1,k_{i+2},\ldots,k_n)$ satisfies $l_i \in
\{k_i,k_i+1\}$. Let us show that the set of penultimate rows with $l_i=k_i$ is
equal to $\mathcal{P}(k_1,\ldots,k_{i-1},k_i,k_i,k_{i+2},\ldots,k_n)$.
It is clear that $l_i=k_i$ implies that the restrictions for
$(l_1,\ldots,l_{i-1})$ are identical for both
$\mathcal{P}(k_1,\ldots,k_{i-1},k_i,k_i+1,k_{i+2},\ldots,k_n)$ and
$\mathcal{P}(k_1,\ldots,k_{i-1},k_i,k_i,k_{i+2},\ldots,k_n)$. For the restrictions of $(l_{i+1},l_{i+2})$
distinguish between $k_i+1 \leq k_{i+2}$, $k_i = k_{i+2}$ and $k_i > k_{i+2}$:

If $k_i+1 \leq k_{i+2}$, then $k_i+1 \leq l_{i+1} \leq k_{i+2}$ on both sides
and the restrictions for $l_{i+2}$ are the same: 
\begin{center}
$
\begin{array}{|c|}
\hline 
\text{Left-hand side of \eqref{oneStepIso}} \\ 
\hline
\begin{array}{cccccccccccccc}
&& k_i &&&& l_{i+1}  \\
& \eqRP && \lRN && \leqRP && \leqRN  \\
k_{i} && < && k_{i}+1 && \leq && k_{i+2} \\ 
\end{array}
\\
\hline\hline
\text{Right-hand side of \eqref{oneStepIso}} \\
\hline
\begin{array}{cccccccccccccc}
&& k_i &&&& l_{i+1}  \\
& \eqRP && \eqRN && \lRP && \leqRN \\
k_{i} && = && k_{i} && \leq && k_{i+2} \\ 
\end{array}
\\
\hline
\end{array}
$
\end{center}

If $k_i = k_{i+2}$, then
$\mathcal{P}(k_1,\ldots,k_{i-1},k_i,k_i,k_{i+2},\ldots,k_n)$ is empty, and
each element of $\mathcal{P}(k_1,\ldots,k_{i-1},k_i,k_i+1,k_{i+2},\ldots,k_n)$
with $l_i=k_i$ would have to satisfy $l_i=l_{i+1}=l_{i+2}=k_i$: 
\begin{center}
$
\begin{array}{|c|}
\hline 
\text{Left-hand side of \eqref{oneStepIso}} \\ 
\hline
\begin{array}{cccccccccccccc}
&& k_i &&&& k_i &&&& k_i  \\
& \eqRP && \lRN && \gRP && \eqRN  \\
k_{i} && < && k_{i}+1 && > && k_{i} \\ 
\end{array}
\\
\hline\hline
\text{Right-hand side of \eqref{oneStepIso}} \\
\hline
\begin{array}{cccccccccccccc}
&& k_i &&&& \lightning &&&&   \\
& \eqRP && \eqRN &&  &&   \\
k_{i} && = && k_{i} && = && k_{i} \\ 
\end{array}
\\
\hline
\end{array}
$
\end{center}
But, since a GMT can not contain three consecutive equal entries, there is
also no element in $\mathcal{P}(k_1,\ldots,k_{i-1},k_i,k_i+1,k_i,\ldots,k_n)$
with $l_i = k_i$.

If $k_i > k_{i+2}$, then $k_i \geq l_{i+1} \geq k_{i+2}$ on both sides and the
restrictions for $l_{i+2}$ are the same:
\begin{center}
$
\begin{array}{|c|}
\hline 
\text{Left-hand side of \eqref{oneStepIso}} \\ 
\hline
\begin{array}{cccccccccccccc}
&& k_i &&&& l_{i+1}  \\
& \eqRP && \lRN && \gRP && \geqRN  \\
k_{i} && < && k_{i}+1 && > && k_{i+2} \\ 
\end{array}
\\
\hline\hline
\text{Right-hand side of \eqref{oneStepIso}} \\
\hline
\begin{array}{cccccccccccccc}
&& k_i &&&& l_{i+1}  \\
& \eqRP && \eqRN && \geqRP && \geqRN \\
k_{i} && = && k_{i} && > && k_{i+2} \\ 
\end{array}
\\
\hline
\end{array}
$
\end{center}
The entry $l_{i+1}$ is involved in a sign-change on both sides (note the
special case $l_{i+1}=k_i$, where $l_{i+1}$ is a newcomer on the left-hand side
and in a sign-changing pair on the right-hand side).

The restrictions for $(l_{i+3},\ldots,l_{n-1})$ are clearly the same for both
sides. Symmetrically, one can also see that the set
$\mathcal{P}(k_1,\ldots,k_{i-1},k_i,k_i+1,k_{i+2},\ldots,k_n)$ restricted to
$l_i = k_{i+1}$ is the same as
$\mathcal{P}(k_1,\ldots,k_{i-1},k_i+1,k_i+1,k_{i+2},\ldots,k_n)$, concluding the
combinatorial proof of \eqref{sasId1} for arbitrary
$(k_1,\ldots,k_n)\in\mathbb{Z}^n$.

A natural question could now be, whether similar identities hold if the
difference between $k_{i+1}$ and $k_i$ is larger. For fixed integers
$k_1,\ldots,k_{i-1},k_{i+2},\ldots,k_n$, let 
\[
t_n(k_i,k_{i+1}):=\alpha(n;k_1,\ldots,k_{i-1},k_i,k_{i+1},k_{i+2},\ldots,k_n).
\]
Similarly - with a bit more patience - one can also show the
identity
\begin{multline}
\label{sasId2}
t_n(k_i,k_i+2)\\=t_n(k_i,k_i)+t_n(k_i+1,k_i+1)+t_n(k_i+2,k_i+2)+t_n(k_i+2,k_i+1)+t_n(k_i+1,k_i)
\end{multline}
combinatorially. Both \eqref{sasId1} and \eqref{sasId2} are special cases of the
following identity: Let $V_{x,y}$
be the operator defined as
\[
V_{x,y}f(x,y):=f(x-1,y)+f(x,y+1)-f(x-1,y+1).
\]
The function $f_i(k_1,\ldots,k_n):=V_{k_i,k_{i+1}} \alpha(n;k_1,\ldots,k_n)$
then satisfies
\begin{equation}
\label{sasId}
f_i(k_1,\ldots,k_n)=-f_i(k_1,\ldots,k_{i-1},k_{i+1}+1,k_i-1,k_{i+2},\ldots,k_n).
\end{equation}

Setting $k_{i+1}=k_i-1$ in \eqref{sasId} immediately implies \eqref{sasId1}.
Equation \eqref{sasId2} is then the special case $k_{i+1}=k_i-2$ in
\eqref{sasId}. A similar shift-antisymmetry property for \emph{Gelfand-Tsetlin
Patterns} (Monotone Triangles without the condition of strict increase along
rows) was shown bijectively in a recent work (\cite{FischerSeqLabTrees}).
It would be interesting to give a bijective proof of \eqref{sasId} in the general
case (an algebraic proof was given in \cite{FischerNumberOfMT}). 

In \cite{FischerRieglerDMT} we showed the surprising identity
\begin{equation}
\label{combRec}
A_n := \alpha(n;1,2,\ldots,n) = \alpha(2n;n,n,n-1,n-1,\ldots,1,1)
\end{equation}
algebraically and gave initial thoughts on how a bijective proof could succeed.
Let us conclude with a list of related identities -- all of them are up to this
point conjectured using mathematical computing software. As Theorem
\ref{gmtTheorem} provides a combinatorial interpretation of these
identities, bijective proofs are of high interest.

\begin{conjecture}[\cite{FischerRieglerDMT}]
Let $n \geq 1$. Then
\begin{equation}
\alpha(2n+1;2n+1,2n,\ldots,1) = (-1)^n \alpha(n;
2,4,\ldots,2n)
\end{equation}
seems to hold, whereby $\alpha(n;2,4,\ldots,2n)$ is known to be the number of
Vertically Symmetric ASMs of size $2n+1$. 
\end{conjecture}

\begin{conjecture}
Let $n \geq 1$. Then
\begin{equation}
\alpha(n;2,4,\ldots,2n) = \alpha(2n;2n,2n,2n-2,2n-2,\ldots,2,2)
\end{equation}
seems to hold.
\end{conjecture}

\begin{conjecture}
Let $n \geq 1$. Then
\begin{align}
  &A_n = \alpha(n+i;1,2,\ldots,i,1,2,\ldots,n), \quad i=0,\ldots,n,\\
  &A_n = (-1)^n \alpha(2n+1;1,2,\ldots,n+1,1,2,\ldots,n)
\end{align}
seems to hold. Furthermore, the numbers
\[
W_{n,i}=\alpha(2n+1;i,2,\ldots,n+1,1,2,\ldots,n), \quad i=1,\ldots,3n+2 
\]
seem to satisfy the symmetry $W_{n,i}=W_{n,3n+3-i}$.
\end{conjecture}

\begin{conjecture}
Let $n \geq 2$. Then
\begin{equation}
\label{oneDescConj}
A_n = \alpha(n+2;1,2,\ldots,i+1,i,i+1,\ldots,n), \quad i=1,\ldots,n-1
\end{equation}
seems to hold.
\end{conjecture}

Further computational experiments led to the conjecture that \eqref{combRec} and
\eqref{oneDescConj} have the following joint generalization:

\begin{conjecture}
Let $n \geq 1$.  Then
\begin{equation}
\label{revDupConj}
\resizebox{.9\hsize}{!}{$
A_n = 
\alpha(n+k;1,\ldots,i-1,i+k-1,i+k-1,i+k-2,i+k-2,\ldots,i,i,i+k,i+k+1,\ldots,n)
$}
\end{equation}
seems to hold for $i=1,\ldots,n-k+1$, $k=1,\ldots,n$.
\end{conjecture}
In words, the last identity takes a subsequence
$(i,i+1,\ldots,i+k-1)$ of length $k$ of $(1,2,\ldots,n)$, reverses
the order, duplicates each entry and puts the subsequence back.
Identity \eqref{combRec} is thus the special case of \eqref{revDupConj}
where $k=n$.
Applying \eqref{sasId1} and the fact that a GMT can not contain three consecutive equal
entries, shows that \eqref{oneDescConj} is the special case of
\eqref{revDupConj} with $k=2$:
\begin{multline*}
\alpha(n+2; 1,2,\ldots,i-1,i,i+1,i,i+1,i+2,\ldots,n) \\
=\alpha(n+2; 1,2,\ldots,i-1,i,i,i,i+1,i+2,\ldots,n)+\alpha(n+2;
1,2,\ldots,i-1,i+1,i+1,i,i+1,i+2,\ldots,n) \\
=\alpha(n+2;1,2,\ldots,i-1,i+1,i+1,i,i,i+2,\ldots,n) +
\alpha(n+2;1,2,\ldots,i-1,i+1,i+1,i+1,i+1,i+2,\ldots,n) \\
=\alpha(n+2;1,2,\ldots,i-1,i+1,i+1,i,i,i+2,\ldots,n).
\end{multline*}

From the correspondence between ASMs of size $n$ and Monotone Triangles with
bottom row $(1,2,\ldots,n)$, it follows that
$\alpha(n-1;1,2,\ldots,i-1,i+1,\ldots,n)$ is equal to the number of ASMs
of size $n$ with the first row's unique $1$ in column $i$ -- denoted $A_{n,i}$.
In the following conjecture we analogously remove the $i$-th argument of the
right-hand side in \eqref{oneDescConj}:

\begin{conjecture}
Let $n \geq 1$. Then
\begin{equation}
\label{holeOneDescConj}
\alpha(n+1; 1,2,\ldots,i-1,i+1,i,i+1,\ldots,n) = -\sum_{j=1}^n (j-i) A_{n,j},
\quad i = 1,\ldots,n-1
\end{equation}
seems to hold.
\end{conjecture}
As a note on how we found \eqref{holeOneDescConj}, let us
prove the case $i=1$: Each penultimate row $(l_1,\ldots,l_n)$ of a GMT with bottom row
$(2,1,2,\ldots,n)$ satisfies $l_1=l_2=1$ by Condition \eqref{gmtDef3} of GMTs.
Taking Conditions \eqref{gmtDef1} and \eqref{gmtDef2} into account, Lemma
\ref{sumOpGMT} implies that
\[
\alpha(n+1;2,1,2,\ldots,n) = -\sum_{p=2}^n
\alpha(n;1,1,2,\ldots,p-1,p+1,\ldots,n).
\]
Each penultimate row $(m_1,\ldots,m_{n-1})$ of a GMT with bottom row
$(1,1,2,\ldots,p-1,p+1,\ldots,n)$ satisfies $m_1=1, m_2=2, \ldots,
m_{p-1}=p-1$. Applying Lemma \ref{sumOpGMT} again yields the claimed equation:
\[
\alpha(n+1;2,1,2,\ldots,n) = -\sum_{p=2}^n \sum_{j=p}^n A_{n,j} \\
= \sum_{j=2}^n (j-1) A_{n,j}.
\]
For general $i$, the set of GMTs with bottom row
$(1,2,\ldots,i-1,i+1,i,i+1,\ldots,n)$ can be written as disjoint union of those
with structure

\begin{center}
\scalebox{0.8}{
$
\begin{array}{ccccccccccccccccccccc}
S_1: && l_1  && \cdots && l_{i-2} && i+1 && i+1 && i && l_{i+2} && \cdots && l_n 
\\
&1  && \cdots && i-2 && i-1 && i+1 && i && i+1 && i+2 && \cdots && n, \\ 
\mbox{}\\
S_2: && l_1  && \cdots && l_{i-2} && i+1 && i && i && l_{i+2} && \cdots && l_n 
\\
&1  && \cdots && i-2 && i-1 && i+1 && i && i+1 && i+2 && \cdots && n, \\ 
\mbox{}\\
S_3: && l_1  && \cdots && l_{i-2} && i-1 && i && i && l_{i+2} && \cdots && l_n 
\\
&1  && \cdots && i-2 && i-1 && i+1 && i && i+1 && i+2 && \cdots && n. \\ 
\end{array}
$
}
\end{center}
Similar to the case $i=1$, one can see that the signed enumeration of
GMTs with structure $S_3$ is equal to 
\[
-\sum_{j=i+1}^n (j-i) A_{n,j}.
\]
Proving that the signed enumeration of GMTs with structure $S_1$ and $S_2$
yields 
\[
-\sum_{j=1}^{i-1} (j-i) A_{n,j}
\]
remains an open problem. The following conjectures are also related to
\eqref{oneDescConj} by removing the $(i-1)$-st argument of the right-hand side.

\begin{conjecture}
Let $n \geq 4$. Then
\begin{equation}
\label{conj13434}
\alpha(n+1;1,3,4,3,4,5,\ldots,n) = \frac{n+4}{2} A_{n-1}
\end{equation}
seems to hold.
\end{conjecture}
As an immediate consequence of Theorem \ref{gmtTheorem} we obtain (the known
fact) that the evaluation of $\alpha(n;k_1,\ldots,k_n)$ at integral values is
integral. From the definition of ASMs it follows that Vertically Symmetric ASMs
only exist for odd size. Therefore, reflection along the vertical symmetry axis
is a fixed-point-free involution on the set of even-sized ASMs. So, the
number of even-sized ASMs is even and the the right-hand side of
\eqref{conj13434} is an integer too.

Using C. Krattenthaler's Mathematica package RATE, we were able to find more
conjectured formulas similar to \eqref{conj13434}:

\begin{conjecture}
\begin{align*}
\alpha(n+1;1,2,4,5,4,5,\ldots,n) = \frac{n^3+7n^2+10n-36}{8n-12} A_{n-1}, \quad
n \geq 5, \\
\alpha(n+1;1,2,3,5,6,5,6,\ldots,n) = \frac{n^4+12n^3++53n^2+54n-288}{48n-72}
A_{n-1}, \quad n \geq 6. \\
\end{align*}
\end{conjecture}

In general, this leads to the following conjecture:
\begin{conjecture}
Let $n \geq k \geq 4$. Then there exist polynomials $p_k(n)$ and $q_k(n)$ with
$\deg p_k - \deg q_k = k-3$ such that
\[
\alpha(n+1;1,2,\ldots,k-3,k-1,k,k-1,k,\ldots,n) = \frac{p_k(n)}{q_k(n)} A_{n-1}.
\]
\end{conjecture}

\bibliographystyle{alpha}
\bibliography{bib120718}

\end{document}